\documentclass[a4paper, 10pt]{article}

\usepackage[T1]{fontenc}
\usepackage[utf8]{inputenc}
\usepackage{lmodern}
\usepackage[UKenglish]{babel}

\usepackage[mathlines, pagewise, modulo]{lineno} 
\usepackage{amsmath, amssymb, amsthm, amsfonts} 
\usepackage{etoolbox} 
\newcommand*\linenomathpatch[1]{%
  \cspreto{#1}{\linenomath}%
  \cspreto{#1*}{\linenomath}%
  \csappto{end#1}{\endlinenomath}%
  \csappto{end#1*}{\endlinenomath}%
}
\newcommand*\linenomathpatchAMS[1]{%
  \cspreto{#1}{\linenomathAMS}%
  \cspreto{#1*}{\linenomathAMS}%
  \csappto{end#1}{\endlinenomath}%
  \csappto{end#1*}{\endlinenomath}%
}

\expandafter\ifx\linenomath\linenomathWithnumbers
  \let\linenomathAMS\linenomathWithnumbers
  \patchcmd\linenomathAMS{\advance\postdisplaypenalty\linenopenalty}{}{}{}
\else
  \let\linenomathAMS\linenomathNonumbers
\fi

\linenomathpatch{equation}
\linenomathpatchAMS{gather}
\linenomathpatchAMS{multline}
\linenomathpatchAMS{align}
\linenomathpatchAMS{alignat}
\linenomathpatchAMS{flalign}

\makeatletter
\patchcmd{\mmeasure@}{\measuring@true}{
  \measuring@true
  \ifnum-\linenopenaltypar>\interdisplaylinepenalty
    \advance\interdisplaylinepenalty-\linenopenalty
  \fi
  }{}{}
\makeatother
\usepackage{changepage} 
\rightlinenumbers

\usepackage{changepage}
\usepackage[inline, shortlabels]{enumitem}
\usepackage{framed}

\usepackage[a4paper, left=4cm, right=4cm]{geometry}

\usepackage[textsize=scriptsize, bordercolor=lightgray, linecolor=lightgray, backgroundcolor=white]{todonotes}

\usepackage{csquotes}
\usepackage[backend=biber, style=alphabetic, doi=false, isbn=false, url=false]{biblatex}
\renewbibmacro{in:}{}
\usepackage{bookmark}
\hypersetup{colorlinks=true, linkcolor=blue, citecolor=blue}

\usepackage{stmaryrd} 
\usepackage{mathabx} 

\theoremstyle{plain}
\newtheorem*{theorem*}{Theorem}
\newtheorem*{conjecture*}{Conjecture}
\newtheorem*{fact*}{Fact}

\newtheorem{question}{Question}
\newtheorem{theorem}{Theorem}
\newtheorem{corollary}{Corollary}[theorem]

\newcommand{\currentstatementname}{}
\newtheorem*{genericclaim}{\currentstatementname}
\newenvironment{citedstatement}[1]{\renewcommand{\currentstatementname}{#1}\begin{genericclaim}}{\end{genericclaim}}

\newtheorem*{definition*}{Definition}
\newtheorem*{notation}{Notation}
\newtheorem{proposition}{Proposition}[theorem]

\theoremstyle{definition}
\newtheorem*{remark}{Remark}
\newtheorem*{remarks}{Remarks}
\newtheorem{step}{Step}[proposition]
\newenvironment{proofclaim}{\begin{proof}}{\end{proof}}

\newcommand{\fa}{\mathfrak{a}}
\newcommand{\fb}{\mathfrak{b}}
\newcommand{\fc}{\mathfrak{c}}
\newcommand{\fd}{\mathfrak{d}}
\newcommand{\fg}{\mathfrak{g}}

\newcommand{\fh}{\mathfrak{h}}
\newcommand{\fr}{\mathfrak{r}}

\newcommand{\fu}{\mathfrak{u}}

\newcommand{\fsl}{\mathfrak{sl}}

\newcommand{\fgl}{\mathfrak{gl}}

\newcommand{\GL}{\mathrm{GL}}
\newcommand{\PGL}{\mathrm{PGL}}
\newcommand{\SL}{\mathrm{SL}}

\newcommand{\bK}{\mathbb{K}}

\newcommand{\bN}{\mathbb{N}}

\newcommand{\ad}{\operatorname{ad}}
\newcommand{\chara}{\operatorname{char}}
\newcommand{\cork}{\operatorname{codim}}
\newcommand{\End}{\operatorname{End}}
\newcommand{\DefEnd}{\operatorname{DefEnd}}
\newcommand{\Id}{\operatorname{Id}}

\newcommand{\rk}{\operatorname{dim}}

\newcommand{\Mod}[1]{#1\text{-}\mathbf{Mod}}
\newcommand{\Vect}[1]{#1\text{-}\mathbf{Vect}}

\title{Soluble Lie rings of finite Morley rank}
\author{Adrien Deloro and Jules Tindzogho Ntsiri}

\begin{document}

\maketitle

\begin{abstract}
We do two things.
\begin{enumerate*}
\item
As a corollary to a stronger linearisation result (Theorem~\ref{t:linearisation}), we prove the finite Morley rank version of the Lie-Kolchin-Malcev theorem on Lie algebras 
(Corollary~\ref{c:Lie}).
\item
We classify Lie ring actions on modules of characteristic $\neq 2, 3$ and Morley rank $2$ (Theorem~\ref{t:module2}).
\end{enumerate*}
\end{abstract}

\section{Introduction}

We keep the introduction short and refer to \cite{DTLie} for more details.
Model-theoretic algebra on Lie rings is young.
Zilber \cite{ZUncountably} had understood that the characteristic~$0$ case reduces to Lie algebra. Nesin \cite{NNonassociative} saw Lie rings as part of a general study of non-associative rings. Then Rosengarten \cite{RAleph} studied with success Lie rings of Morley rank $1, 2, 3$. Our \cite{DTLie} proved that in characteristic $\neq 2, 3$, there are no simple Lie rings of Morley rank $4$.

In the present paper we solve two basic problems, which should become pervasive tools in the theory of Lie rings of finite Morley rank.
\begin{enumerate}
\item
The analogue of Lie's theorem, viz.~the fact that the commutator algebra of a complex, finite-dimensional, soluble Lie algebra is nilpotent.
This first step towards \cite[Question~10]{DTLie} is taken in \S~\ref{S:linearisation}.
\item
Classification of $2$-dimensional modules for abstract Lie rings of finite Morley rank.
This solution to \cite[Question~8]{DTLie} is given in \S~\ref{S:module2}.
\end{enumerate}

\subsection{Present results (for quick reference)}

There is terminology and notation in \S~\ref{s:terminology}.

\begin{citedstatement}{Theorem~\ref{t:linearisation} (\S~\ref{s:linearisation:main})}
Let $\fg$ be a Lie ring of finite Morley rank (not necessarily connected). Let $V$ be a definable, irreducible, faithful $\fg$-module. 
Suppose $\chara V > \rk V$.
Suppose that $\fg$ has an infinite abelian ideal $\fa \trianglelefteq \fg$.
Then $\fa \sqsubseteq Z(\fg)$ and the configuration is definably linear: there is an infinite field $\bK$ such that $V \in \Vect{\bK}$, $\fg \hookrightarrow \fgl(V\colon \Vect{\bK})$, and $\fa \hookrightarrow \bK \Id_V$, all definably.
\end{citedstatement}

\begin{citedstatement}{Corollary~\ref{c:solubleaction} (\S~\ref{s:linearisation:Lie})}
Let $\fr$ be a connected, soluble Lie ring of finite Morley rank. Let $V$ be a definable, irreducible, faithful $\fr$-module. Suppose $\chara V > \rk V$. Then $\fr$ is abelian and the configuration is definably linear: there is an infinite field $\bK$ such that $V\simeq \bK_+$ and $\fr \hookrightarrow \bK \Id_V$, all definably.
\end{citedstatement}

\begin{citedstatement}{Corollary~\ref{c:Lie} (\S~\ref{s:linearisation:Lie})}
Let $\fr$ be a connected, soluble Lie ring of finite Morley rank. Suppose $(\fr ; +)$ has prime exponent $> \rk \fr$. Then:
\begin{enumerate}[label=(\roman*)]
\item
there is a largest definable, connected, nilpotent ideal $F^\circ(\fr) \trianglelefteq \fr$;
\item
$\fr'$ is nilpotent, viz.~$\fr' \sqsubseteq F^\circ(\fr)$;
\item
if $\fu \sqsubseteq \fr$ is a definable, connected subring consisting of $\fr$-nilpotent elements (viz.~$(\forall a \in \fu)(\exists n \in \bN)(\ad_a^n(\fr) = 0)$), then $\fu \sqsubseteq F^\circ(\fr)$.
\end{enumerate}
\end{citedstatement}

\begin{citedstatement}{Theorem~\ref{t:module2} (\S~\ref{S:module2})}
Let $\fg$ be a connected Lie ring of finite Morley rank acting definably and faithfully on an irreducible module $V$ of characteristic $\neq 2, 3$.
Suppose $\rk V = 2$. Then there is a definable field $\bK$ such that, all definably:
\begin{itemize}
\item
either $V \simeq \bK_+$ and $\fg \hookrightarrow \bK\Id_V$;
\item
or $\fg \simeq \fsl_2(\bK)$ in its natural action on $V\simeq \bK^2$;
\item
or $\fg \simeq \fgl_2(\bK)$ in its natural action on $V\simeq \bK^2$.
\end{itemize}
\end{citedstatement}

\subsection{Terminology and notation}\label{s:terminology}

\begin{itemize}
\item
Inside Lie rings, we use $\leq$ for subgroups, $\sqsubseteq$ for subrings, and $\trianglelefteq$ for ideals.
\item
For $V$ an abelian group, $(\End(V) ; +, \cdot)$ is an associative ring.
When there is a notion of definability around, we let $\DefEnd(V)$ be the subring of \emph{definable} endomorphisms; $\DefEnd(V)$ itself need not be definable.
\item
Both $\End(V)$ and $\DefEnd(V)$ bear a Lie ring structure by letting $\llbracket f, g\rrbracket = f\circ g - g \circ f$.
\item
A Lie ring $\fg$ \emph{acts} on $V$ if there is a morphism $\rho\colon (\fg ; +, [\cdot, \cdot]) \to (\End(V) ; +, \llbracket\cdot, \cdot\rrbracket)$. When there is a notion of definability around, the action is \emph{definable} if $(\fg ; +, [\cdot, \cdot])$, $(V ; +)$, and the map $\fg \times V \to V$ taking $(g, v)$ to $\rho(g)(v)$ are definable. We omit $\rho$ from notation.
\item
A Lie ring acting definably on $V$ generates an \emph{invariant} associative subring of $\DefEnd(V)$. Here, \emph{invariance} is a model-theoretic notion generalising type-definability. For our purpose, $\bigvee$-definability is a decent proxy.
\item
In case $V$ is definable and connected we call $V$ a \emph{$\fg$-module}; we avoid the phrase otherwise. In particular, `$V$ is an irreducible $\fg$-module' means that $V$ has no definable, connected, non-trivial, proper, $\fg$-invariant subgroup.
\item
Let $V$ be a definable, irreducible $\fg$-module in a theory of finite Morley rank. Then:
\begin{itemize}
\item
either $(V ; +)$ has prime exponent $p > 0$, called the \emph{characteristic} of $V$;
\item
or $(V ; +)$ is divisible, in which case we say $V$ has \emph{characteristic $0$}.
\end{itemize}
This follows from Macintyre's theorem on abelian groups of ordinal Morley rank. (By rigidity phenomena, there can be no non-trivial definable module structure on a divisible torsion group; also, by \cite[Theorem~0]{DTLie}, we have no remaining interest in characteristic $0$.)
\item
$\Vect{\bK}$ is the category of $\bK$-vector spaces.
\item
The general linearisation result \cite{DZilber} has a consequence on modules of Morley rank $1$ \cite[Corollary~1]{DTLie}: they linearise without solubility.
\item
$\rk$ stands for Morley rank. It has been noted that a weaker assumption than finiteness of the Morley rank might suffice for model-theoretic algebra on Lie rings \cite[Question~13]{DTLie}. If would be interesting (though possibly very challenging) to see what happens if one allows definable field derivations, as opposed to the finite Morley rank context \cite[Lemma~G]{DTLie}.
\end{itemize}

\section{Linearisation, and Lie's theorem}\label{S:linearisation}

Our ad hoc \cite[Corollary~2]{DTLie} was a disappointingly special case of a desirable result, viz.~`Lie's theorem' for abstract Lie rings of finite Morley rank.
Its group analogue, sometimes called `Lie-Kolchin-Malcev theorem', states that the commutator subgroup of a connected, soluble group of finite Morley rank is nilpotent. This was first proved for groups of finite Morley rank by Nesin \cite{NSolvable}.
We obtain the Lie ring version (Corollary~\ref{c:Lie} below) as a consequence of Theorem~\ref{t:linearisation}, a linearisation result similar to \cite[Theorem~3.9]{PStable}. Poizat had extracted the latter from \cite{NSolvable}.

Though new, Corollary~\ref{c:Lie} is not unexpected. It is `presumed' without a proof as \cite[Theorem~18]{NNonassociative}. However Nesin did not make the necessary assumptions on the characteristic (personal conversation confirmed he did not know about low-characteritic counterexamples).
Also, in our context, lack of conjugacy prevents us from retrieving a global linear structure; this sharply contrasts with \cite{NSolvable} and \cite[Theorem~3.9]{PStable}.
Thus Theorem~\ref{t:linearisation} is not as straightforward as it looks.

\subsection{Linearisation theorem}\label{s:linearisation:main}

\begin{theorem}\label{t:linearisation}
Let $\fg$ be a Lie ring of finite Morley rank (not necessarily connected). Let $V$ be a definable, irreducible, faithful $\fg$-module. 
Suppose $\chara V > \rk V$.
Suppose that $\fg$ has an infinite abelian ideal $\fa \trianglelefteq \fg$.
Then $\fa \sqsubseteq Z(\fg)$ and the configuration is definably linear: there is an infinite field $\bK$ such that $V \in \Vect{\bK}$, $\fg \hookrightarrow \fgl(V\colon \Vect{\bK})$, and $\fa \hookrightarrow \bK \Id_V$, all definably.
\end{theorem}
\begin{proof}
Let $\hat\fa = Z(C_\fg(\fa)) \sqsupseteq \fa$. Then $\hat\fa \trianglelefteq \fg$ is an infinite abelian ideal, and is definable. So we may assume that $\fa$ is definable.
Both $\fg$ and $\fa$ generate invariant, unbounded, associative subrings of $\DefEnd(V)$.
By \cite[Corollary~1]{DZilber}, it is enough to prove $\fa \sqsubseteq Z(\fg)$. So suppose not. We seek a contradiction.

\begin{notation}
Let $W \leq V$ be $\fa$-irreducible.
\end{notation}

\begin{step}\label{st:Wlinear}\leavevmode
\begin{enumerate}[label={\itshape (\roman*)},series=claims]
\item\label{st:W:i:noideal}
No non-zero ideal of $\fg$ may centralise $W$.
\item\label{st:W:i:h}
There is $h \in \fg$ such that $[h, \fa]$ does not centralise $W$.
\item\label{st:W:i:W}
$W$ has no non-trivial, finite $\fa$-invariant subgroups.
\end{enumerate}
\end{step}
\begin{proofclaim}
\begin{enumerate}[label={\itshape (\roman*)},series=proofs]
\item
Let $0 < \fb \trianglelefteq \fg$ be any non-zero ideal. Then $C_V^\circ(\fb)$ is an $\fg$-submodule of $V$. By faithfulness, it is proper, so by $\fg$-irreducibility, it is trivial. Hence $\fb$ may not centralise $W$.
\item
Let $\fb = [\fg, \fa] \trianglelefteq \fg$.
We supposed $\fa \not\sqsubseteq Z(\fg)$, so $\fb \neq 0$.
By~\ref{st:W:i:noideal}, $\fb$ acts non-trivially on $W$. So there is $h \in \fg$ such that $[h, \fa]$ does not centralise $W$.
\item
Let $\fb = \fa^\circ \trianglelefteq \fg$. Since $\fa$ is infinite, $\fb \neq 0$. By~\ref{st:W:i:noideal}, $\fb$ acts non-trivially on $W$, so $\overline{\fa} = \fa/C_\fa(W)$ is infinite. Let $\bK = C_{\DefEnd(W)}(\overline{\fa})$. By abelianity, $\overline{\fa} \leq \bK$.
Linearising \cite{DZilber}, $\bK$ is an infinite definable field.
So $\overline{\fa}$ acts by scalars, hence freely, on $W$. As it is infinite, there can be no finite non-trivial $\fa$-invariant subgroup of $W$.
\qedhere
\end{enumerate}
\end{proofclaim}

An element $h$ like in~\ref{st:W:i:h} is now fixed.
For $a \in \fa$, let
$a^{[i]} = \ad_h^i(a)$. In particular, $a' = [h, a]$.
In this notation, for $i \geq 0$ one has in $\DefEnd(V)$:
\[h^i a = \sum_{j = 0}^i \binom{i}{j} a^{[j]}h^{i-j}.\]

\begin{notation}
Set $S_0 = W$ and for $i \geq 1$, let $S_i = \sum_{j = 0}^{i-1} h^j W$.
\end{notation}

The sum giving $S_i$ thus has $i$ terms. (For the computations below, this notation is slightly preferable.)

\begin{step}\leavevmode
\begin{enumerate}[resume*=claims]
\item\label{st:S:i:Si}
Each $S_i$ is an $\fa$-module.
\item\label{st:S:i:hi}
If $w \in W\setminus \{0\}$ and $i \geq 0$ are such that $h^i w \in S_i$, then $S_{i+1} = S_i$.
\item\label{st:S:i:k}
There is $k \leq \rk V$ such that $\sum_{n \geq 0} h^n W = S_k = \bigoplus_{i = 0}^{k-1} h^i W > S_{k-1}$.
\end{enumerate}
\end{step}
\begin{proofclaim}\leavevmode
\begin{enumerate}[resume*=proofs]
\item
Definability and connectedness are obvious; $\fa$-invariance is by induction. Indeed, let $i \geq 0$ be such that $S_i$ is $\fa$-invariant. Let $a \in \fa$ and $w \in W$. One has:
\[a h^i \cdot w = h^i a \cdot w - \sum_{j = 1}^i \binom{i}{j} a^{[j]} h^{i - j} \cdot w \in S_{i+1},\]
which proves $\fa$-invariance of $S_{i+1}$.
\item
Let $W' = \{w \in W: h^i w \in S_i\} \neq 0$. 
For $w \in W'$ and
$a \in \fa$ one has:
\[h^i a \cdot w = \sum_{j = 0}^i \binom{i}{j} a^{[j]}h^{i-j} \cdot w.\]
Each $S_{i-j}$ is an $\fa$-module by~\ref{st:S:i:Si}, so $h^i a \cdot w \in S_i$. Hence $W'$ is a non-zero, definable, $\fa$-invariant subgroup of $W$. By~\ref{st:W:i:W}, it is infinite; so $\fa$-irreducibility of $W$ implies $(W')^\circ = W$. Hence $h^i W \leq S_i$.
\item
If the sum defining $S_i$ is direct, then $i \leq \rk V$.
So there is $k$ maximal such that the sum giving $S_k$ is direct. Then $k \leq \rk V$. Moreover, $S_k \cap h^k W \neq 0$, so there is $w \in W\setminus \{0\}$ with $h^k w \in S_k$. Then $S_{k+1} = S_k$ by~\ref{st:S:i:hi}, and $S_\infty = S_k$.
\qedhere\end{enumerate}
\end{proofclaim}

\begin{notation}\leavevmode
\begin{itemize}
\item
With $k$ as in~\ref{st:S:i:k}, let $X = S_k/S_{k-1}$.
\item
Let $\bK = C_{\DefEnd(X)}(\fa)$.
\end{itemize}
\end{notation}


\begin{step}\leavevmode
\begin{enumerate}[resume*=claims]
\item\label{st:X:i:XW}
Taking $w$ to $(h^{k-1} \cdot w \mod S_{k-1})$ defines an isomorphism $W \simeq X$ of $\fa$-modules.
\item\label{st:X:i:K}
$\bK$ is an infinite, definable field.
\item\label{st:X:i:eta}
$h$ induces a natural morphism $\eta \in \DefEnd(X)$; for $a \in \fa$, one has $\eta a = a \eta + k a'$.
\end{enumerate}
\end{step}
\begin{proofclaim}\leavevmode
\begin{enumerate}[resume*=proofs]
\item
Consider the additive morphism $W \to X$ taking $w$ to $(h^{k-1} \cdot w \mod S_{k-1})$. Since $S_k = h^{k-1} W + S_{k-1}$, it is onto.
Its $\fa$-covariance results from $h^{k-1} a = \sum_{i = 0}^{k-1} \binom{k-1}{i} a^{[i]} h^{k-1-i}$.
If $w \in W \setminus\{0\}$ is such that $h^{k-1} w \in S_{k-1}$, then~\ref{st:S:i:hi} implies $S_{k-1} = S_k$, against the definition of $k$.
So the morphism is also injective.
\item
By~\ref{st:X:i:XW}, $W \simeq X$ as $\fa$-modules. By construction, $X$ is therefore $\fa$-irreducible. The image $\fa/C_\fa(X)$ is infinite by~\ref{st:W:i:noideal}. Linearising the abelian action, $\bK$ is an infinite, definable field.
\item
Let $x \in X$, say $x = (h^{k-1} w \mod S_{k-1})$ with $w \in W$. We let $\eta(x) = (h^k w \mod S_{k-1})$. This is well-defined. Indeed, if $h^{k-1} w \in S_{k-1}$, then~\ref{st:S:i:hi} implies $w = 0$, so $h^k w = 0$ and $(h^k \mod S_{k-1}) = 0$.

Let $x \in X$, say $x = (h^{k-1} w \mod S_{k-1})$ with $w \in W$. Then $h^{k-1} a w \equiv a h^{k-1} w~[S_{k-1}]$, so $h^{k-1} a w$ represents $a \cdot x$ and $\eta(a \cdot x) = (h^k a w \mod S_{k-1})$.
Now:
\[h^k a w \equiv a h^k w + k a' h^{k-1} w~[S_{k-1}],\]
so $\eta (a \cdot x) = (ah^k w + k a' h^{k-1} w \mod S_{k-1}) = a \cdot \eta(x) + k a' \cdot x$.
This proves $\eta a = a \eta + k a'$.
\qedhere
\end{enumerate}
\end{proofclaim}

\begin{step}\leavevmode
\begin{enumerate}[resume*=claims]
\item\label{st:X:i:delta}
For $\lambda \in \bK$, let $\delta(\lambda) = \llbracket \eta, \lambda\rrbracket$.
Then $\delta$ is a definable derivation of $\bK$.
\item
Contradiction.
\end{enumerate}
\end{step}
\begin{proofclaim}\leavevmode
\begin{enumerate}[resume*=proofs]
\item
The $\llbracket\cdot, \cdot\rrbracket$-bracket is computed in $\DefEnd(X)$ and therefore makes sense. It yields $\ad_\eta(\lambda)$, which results in a derivation on $\DefEnd(X)$.
One must however check $\delta(\lambda) \in \bK$, viz.~that $\delta(\lambda)$ commutes with the action of $\fa$.

By~\ref{st:X:i:eta}, one has for $\lambda \in \bK = C_{\DefEnd(X)}(\fa)$:
\begin{align*}
(\eta \lambda - \lambda \eta)a & = (\eta a) \lambda - \lambda (\eta a)\\
& = a \eta\lambda  + k a' \lambda - \lambda a \eta - \lambda k a'\\
& = a (\eta \lambda - \lambda \eta).
\end{align*}
Therefore $\delta(\lambda) \in \bK$.
\item
An infinite field of finite Morley rank has no non-trivial definable derivations by \cite[Lemma~G]{DTLie}.
So for $\lambda \in \bK$ one has $\delta(\lambda) = 0$. In particular, for (the image in $\bK$ of) $a \in \fa$, one has $\llbracket\eta, a\rrbracket = 0$. By~\ref{st:X:i:eta}, $k a'$ acts trivially on $X$. By~\ref{st:S:i:k} and by assumption, $k \leq \rk V < \chara V$. So $a'$ acts trivially on $X$, but $X \simeq W~[\Mod{\fa}]$ by~\ref{st:X:i:XW}.
So $[h, \fa]$ acts trivially on $W$, against~\ref{st:W:i:h}.
\qedhere
\end{enumerate}
\end{proofclaim}

This contradiction completes the proof: $\fa \sqsubseteq Z(\fg)$, so we may linearise.
\end{proof}

\subsection{Corollaries}\label{s:linearisation:Lie}

\begin{corollary}\label{c:solubleaction}
Let $\fr$ be a connected, soluble Lie ring of finite Morley rank. Let $V$ be a definable, irreducible, faithful $\fr$-module. Suppose $\chara V > \rk V$. Then $\fr$ is abelian and the configuration is definably linear: there is an infinite field $\bK$ such that $V\simeq \bK_+$ and $\fr \hookrightarrow \bK \Id_V$, all definably.
\end{corollary}
\begin{proof}
Let $\fa$ be the last non-trivial commutator subring. By the theorem, the configuration is linear, say $V \simeq \bK_+^n$ and $\fr \hookrightarrow \fgl_n(\bK)$. By Lie's theorem, $\fr$ (or the Lie algebra it generates, viz.~its $\bK$-linear span) being soluble has a common eigenvector in $V$. Irreducibility forces $n = 1$; so $\fr$ is abelian.
\end{proof}

In the following, one must assume existence of a characteristic. In practice, all our soluble Lie rings will be definable subrings of simple Lie rings of finite Morley rank, so this requirement will be met.

\begin{corollary}\label{c:Lie}
Let $\fr$ be a connected, soluble Lie ring of finite Morley rank. Suppose $(\fr ; +)$ has prime exponent $> \rk \fr$.
Then:
\begin{enumerate}[label=(\roman*)]
\item\label{c:Li:e:i:Fitting}
there is a largest definable, connected, nilpotent ideal $F^\circ(\fr) \trianglelefteq \fr$;
\item
$\fr'$ is nilpotent, viz.~$\fr' \sqsubseteq F^\circ(\fr)$;
\item
if $\fu \sqsubseteq \fr$ is a definable, connected subring consisting of $\fr$-nilpotent elements (viz.~$(\forall a \in \fu)(\exists n \in \bN)(\ad_a^n(\fr) = 0)$), then $\fu \sqsubseteq F^\circ(\fr)$.
\end{enumerate}
\end{corollary}

\begin{proof}
This does follow the group case, as planned by Nesin. Mostly induction on $\rk \fr$.
Good references are \cite[Corollary~9.9 and Theorem~9.21]{BNGroups}.
\end{proof}

\begin{remarks}\leavevmode
\begin{itemize}
\item
As is well-known (and which makes the claim in \cite[Theorem~18]{NNonassociative} look hasty), there are algebraic counterexamples to nilpotence of $\fr'$ if $(\fr ; +)$ has exponent $\leq \rk \fr$.
\item
We introduce only a connected version of the Fitting ideal, so~\ref{c:Li:e:i:Fitting} is obvious. One could wonder what happens with the non-connected version \cite[\S~7.2]{BNGroups}.
\item
There remains of course to develop a theory of Cartan subrings \cite[Question~10]{DTLie}.
Only then should one attack the analogue of the Cherlin-Zilber conjecture. Fortunately Theorem~\ref{t:module2} is too small to require advanced tools.
\end{itemize}
\end{remarks}

\section{Modules of Morley rank 2}\label{S:module2}

The following is an analogue of \cite[Theorem~A]{DActions}. Its proof may look easy; it however relies on \cite{DTLie}, which has no analogue in the group case.

\begin{theorem}\label{t:module2}
Let $\fg$ be a connected Lie ring of finite Morley rank acting definably and faithfully on an irreducible module $V$ of characteristic $\neq 2, 3$.
Suppose $\rk V = 2$. Then there is a definable field $\bK$ such that, all definably:
\begin{itemize}
\item
either $V \simeq \bK_+$ and $\fg \hookrightarrow \bK\Id_V$;
\item
or $\fg \simeq \fsl_2(\bK)$ in its natural action on $V\simeq \bK^2$;
\item
or $\fg \simeq \fgl_2(\bK)$ in its natural action on $V\simeq \bK^2$.
\end{itemize}
\end{theorem}

\subsection{Proof of Theorem~\ref{t:module2}}

\begin{proof}
Let $\fg$ be a connected Lie ring of finite Morley rank. Let $V$ be a definable, faithful, irreducible $\fg$-module of characteristic $\neq 2, 3$ with $\rk V = 2$.

\begin{step}[reductions]\leavevmode
\begin{enumerate}[resume*=claims, start=1]
\item\label{t:module2:i:nonsoluble}
We may suppose that $\fg$ is non-soluble.
\item\label{t:module2:i:rightkernel}
We may suppose that for every $v \in V\setminus \{0\}$, the centraliser $C_\fg^\circ(v)$ is proper.
\item\label{t:module2:i:Cvsoluble}
We may suppose that for every $v \in V\setminus \{0\}$, the centraliser $C_\fg^\circ(v)$ is soluble.
\item\label{t:module2:i:soluble}
Let $\fh \sqsubseteq \fg$ be a definable, connected subring. Suppose there is $v\in V\setminus \{0\}$ such that $\fh \cap C_\fg(v) \neq 0$. Then $\fh$ is soluble iff $V$ is $\fh$-reducible.
\end{enumerate}
\end{step}
\begin{proofclaim}\leavevmode
\begin{enumerate}[resume*=proofs, start=1]
\item
If $\fg$ is soluble, then Corollary~\ref{c:solubleaction} gives the desired description.
\item
Let $R(V) = \{v \in V: \fg \cdot v = 0\}$, a finite subgroup by $\fg$-irreducibility.
Let $\overline{V} = V/R(V)$ and $\pi\colon V \twoheadrightarrow \overline{V}$ be the natural projection.
Let $R_2(V) = \pi^{-1}(R(\overline{V}))$.
Then $R_2(V)$ is finite and $\fg$-invariant, so by connectedness $\fg \cdot R_2(V) = 0$ and $R_2(V) \leq R(V)$. Hence $R(\overline{V}) = 0$.

Suppose the problem is solved for modules with $R(V) = 0$. Then there is a definable field $\bK$ such that $(\fg, \overline{V})$ is either $(\fgl_2(\bK), \bK^2)$ or $(\fsl_2(\bK), \bK^2)$, both in the natural action. Let $v \in V \setminus R(V)$, with image $\overline{v} = \pi(v) \in \overline{V}$. In the natural action, $\fg \cdot \overline{v} = \overline{V}$, so by connectedness $\fg \cdot v = V$. Let $w \in R(V)$. Then there is $a \in \fg$ such that $w = a \cdot v$. So $a \in C_\fg(\overline{v})$. But $C_\fg(v) \sqsubseteq C_\fg(\overline{v})$ and $C_\fg(\overline{v})$ takes $v$ to $R(V)$, so $C_\fg^\circ(\overline{v}) \sqsubseteq C_\fg^\circ(v)$. In the natural action, centralisers of vectors of $\overline{V}$ in $\fg$ are connected.
Hence:
\[a \in C_\fg(\overline{v}) = C_\fg^\circ(\overline{v}) = C_\fg^\circ(v) \sqsubseteq C_\fg(v),\]
meaning $w = a \cdot v = 0$, as wanted.
We may therefore assume $R(V) = 0$.
\item
Induction on $\rk \fg$.
Let $v \in V\setminus \{0\}$ and $\fc = C_\fg^\circ(v)$.
By~\ref{t:module2:i:rightkernel}, $\fc$ is proper.
If non-soluble, then $V$ is $\fc$-irreducible, since otherwise one may linearise in dimension $1$ \cite[Corollary~1]{DTLie} twice and find $\fc'' = 0$ (more details in the proof of~\ref{t:module2:i:soluble}). By induction $\fc$ is either $\fsl_2(\bK)$ or $\fgl_2(\bK)$ in their natural action: a contradiction to $\fc = C_\fg^\circ(v)$.
\item
Suppose $\fh$ is soluble. If $V$ is $\fh$-irreducible, then by Corollary~\ref{c:solubleaction}, $\fh$ acts freely on $V$, against its meeting a centraliser.
Suppose $V$ is $\fh$-reducible, say $0 < W < V$ is an $\fh$-series. Linearising in dimension $1$ twice, $\fh'$ acts trivially on both $W$ and $V/W$, so $\fh'' = 0$. Hence $\fh$ is soluble.
\qedhere
\end{enumerate}
\end{proofclaim}

\begin{step}[real-life cases]\leavevmode
\begin{enumerate}[resume*=claims]
\item\label{t:module2:i:3}
If $\rk \fg = 3$, then we are done.
\item\label{t:module2:i:4}
If $\rk \fg = 4$, then we are done.
\end{enumerate}
\end{step}
\begin{proofclaim}\leavevmode
\begin{enumerate}[resume*=proofs]
\item
Suppose $\rk \fg = 3$. It is non-soluble by~\ref{t:module2:i:nonsoluble}, so by \cite{RAleph} or \cite{DTLie}, $\fg$ has a finite centre and $\fg/Z(\fg) \simeq \fsl_2(\bK)$ for some definable field $\bK$.

We prove $Z(\fg) = 0$ and $\fg \simeq \fsl_2(\bK)$.
Let $\overline{\fg} = \fg/Z(\fg)$.
Let $\overline{h} \in \overline{\fg}$ induce a weight decomposition $\overline{\fg} = E_{-2}(\overline{h}) \oplus E_0(\overline{h}) \oplus E_2(\overline{h})$.
Here $E_k(\overline{h}) = \ker^\circ(\ad_{\overline{h}} - k \Id) \leq \overline{\fg}$ (see \cite[\S~2.1]{DTLie}).
Let $h \in \fg$ lie above $\overline{h}$.
We turn to eigenspaces for $h$ in $\fg$, simply denoted by $E_k$.
Lifting eigenspaces \cite[Lemma~D]{DTLie} and using connectedness, $\fg = E_{-2} \oplus E_0 \oplus E_2$. If $z \in Z(\fg)$, then $z$ decomposes as $z = z_{-2} + z_0 + z_2$; applying $h$ shows $z_{-2} = z_2 = 0$. Hence $Z(\fg) \sqsubseteq E_0$. Fix $e \in E_2 \setminus\{0\}$. Then $e$ is not central, so by connectedness $[E_{-2}, e] = E_0$. In particular, there is $f \in E_{-2}$ such that $[f, e] = z$. But then $[\overline{f}, \overline{e}] = 0$, forcing one to be trivial in $\overline{\fg}$. Since $Z(\fg) \sqsubseteq E_0$, one of $e$ or $f$ is trivial in $\fg$, and so is $z$. Thus $Z(\fg) = 0$ and $\fg \simeq \fsl_2(\bK)$.

Let $\fb \sqsubset \fsl_2(\bK)$ be a Borel subring. If $V$ is $\fb$-irreducible, then by Corollary~\ref{c:solubleaction}, $\fb$ is abelian: a contradiction. So there is a $\fb$-series $0 < W < V$. By Corollary~\ref{c:solubleaction} again, $\fb'$ acts trivially on $W$ and $V/W$. In particular, $\fb'$ acts quadratically, viz.~for $x \in \fb'$ one has $x^2 \cdot V = 0$.

Finally, $V$ has no non-trivial, finite $\fg$-invariant subgroup by~\ref{t:module2:i:rightkernel}.
This is enough to follow the proof of \cite[Variation~11]{DQuadratic}, which definably provides a $\bK$-linear structure on $V$ turning it into the natural representation.
\item
Suppose $\rk \fg = 4$.
By \cite{DTLie}, $\fg$ cannot be simple; so there is an infinite, proper, definable, connected $I \triangleleft \fg$.

If $\rk I = 1$, then $I$ is abelian and we linearise using Theorem~\ref{t:linearisation}: there is a definable field $\bK$ such that $V \in \Vect{\bK}$ and $\fg \hookrightarrow \fgl(V\colon \Vect{\bK})$. Clearly $\rk \bK = 1$ and the $\bK$-linear dimension of $V$ is $2$; we are done.

If $\rk I = 2$, then by \cite[Theorem~4.2.1]{RAleph} $I$, $\fg/I$, and $\fg$ are soluble; against~\ref{t:module2:i:soluble}.

Therefore $\rk I = 3$, and by~\ref{t:module2:i:3} $I \simeq \fsl_2(\bK)$ in its natural action on $V \simeq \bK^2$.

Let $W$ be a $\bK$-vector line in $V \simeq \bK^2$.
Fix $w \in W\setminus \{0\}$. Let $\fc = C_\fg^\circ(w)$, then $\fu = C_I^\circ(w)$ and $\fb = N_I^\circ(\fu)$. Observing the natural representation, $\fu$ is a nilpotent $\bK$-subalgebra of $I$ and $\fb$ the Borel $\bK$-subalgebra of $I$ containing $\fu$. Moreover, $\fu = C_I^\circ(W)$ and $W = C_V^\circ(\fu)$. So $\fc$ normalises $\fu$ and $W$.

If $\rk \fc = 3$, then $\rk \fu = \rk(\fc \cap I)^\circ \geq 2$, a contradiction. Hence $\rk \fc = 2$.
Let $\fd = \fb + \fc$.
Since $\fc \not\sqsubseteq I$, one has $\rk \fd \geq 3$. Since $\fd$ normalises $W$, one has $\rk \fd = 3$ by irreducibility of $\fg$.

Now start from three distinct $\bK$-vector lines in $V$, say $W_1, W_2, W_3$. Fix $w_i \in W_i \setminus \{0\}$ in each and define $\fc_i, \fu_i, \fb_i, \fd_i$ accordingly. For $i \neq j$, let $\fd_{i, j} = (\fd_i \cap \fd_j)^\circ$, which normalises $W_i$ and $W_j$. By Lie's Theorem, $\fd_{i, j}' = 0$, so $\fd_{i, j}$ is abelian. Since $\rk (\fd_{i, j} \cap I) \leq 1$, one has $\rk \fd_{i, j} = 2$, and $(\fd_{i, j}\cap I)^\circ$ is a Cartan $\bK$-subalgebra of $I$.
In particular, $\fd_{i, j}$ does not normalise $W_k$, implying $\fd_{i, j} \not\sqsubseteq \fd_{i, k}$.

Finally let $\fa = (\fd_1 \cap \fd_2 \cap \fd_3)^\circ$; hence $\rk \fa = 1$.
Let $\fh = C_\fg^\circ(\fa)$. Then $\fh$ contains $\fd_{i, j}$ for each pair. So $\fh \cap I$ is a subring of $I$ containing two distinct Cartan subalgebras of $I$: hence $I \leq \fh$. By abelianity of $\fa$ \cite[Theorem~4.1.1]{RAleph}, also $\fa \leq \fh$. Thus $\fa$ is an abelian ideal of $\fg$, and we are back to the first case above.
\qedhere
\end{enumerate}
\end{proofclaim}

\begin{definition*}
A definable, connected subring $\fu \sqsubseteq \fg$ is \emph{$V$-nilpotent} if for all $a \in \fu$, there is an integer $n$ such that $a^n = 0~[\End(V)]$.
\end{definition*}

In this terminology, $\ad$-nilpotence becomes $\fg$-nilpotence, where $\fg$ is the adjoint module.
In general, if $\fg$ is faithful on $V$ and $\fu \sqsubseteq \fg$ is $V$-nilpotent with lengths bounded by $n$, then $\fu$ is $\fg$-nilpotent with lengths bounded by $2n-1$.
(This relies on the formula $\ad_a^k(x) = \sum_{i = 0}^k \binom{k}{i} a^{k-i} x a^i~[\End(V)]$.) We shall use $\fg$-nilpotence of $V$-nilpotent subrings in~\ref{t:module2:i:u1u2} below.

The following manages to elude an $N_\circ^\circ$-theory, which would be still interesting to have, with the proviso that Witt's algebras are simple and $N_\circ^\circ$; see \cite[Question~11]{DTLie}.

\begin{step}[nilpotent analysis]\leavevmode
\begin{enumerate}[resume*=claims]
\item\label{t:module2:i:u}
If $v \in V\setminus\{0\}$, then $C_\fg^\circ(v)$ contains a non-trivial, $V$-nilpotent subring.

Let $\fu$ be a maximal $V$-nilpotent subring.
\item\label{t:module2:i:Wu}
There exists a unique non-trivial, proper $\fu$-submodule $W_\fu < V$. One has $\fu = C_\fg^\circ(W_\fu, V/W_\fu)$ and $W_\fu = C_V^\circ(\fu)$.
\item\label{t:module2:i:bu}
$N_\fg^\circ(\fu) = N_\fg^\circ(W_\fu)$ is a Borel subring of $\fg$.
\end{enumerate}
\end{step}
\begin{proofclaim}\leavevmode
\begin{enumerate}[resume*=proofs]
\item
Let $\fc = C_\fg^\circ(v)$. By~\ref{t:module2:i:Cvsoluble}, $\fc$ is soluble, so by~\ref{t:module2:i:soluble}, it is reducible.
Let $0 < W < V$ be a $\fc$-series, so $\rk W = 1$. Let $\fc_1 = C_\fc(W)$ and $\fc_2 = C_\fc(V/W)$. (We do not use connected components here.)
By definition, $\fu = (\fc_1 \cap \fc_2)^\circ$ is $V$-nilpotent.

We claim $\fc = \fc_1 \cup \fc_2$. Otherwise, linearising in dimension $1$, any $a \in \fc\setminus (\fc_1 \cup \fc_2)$ acts like a non-zero scalar (hence freely) on both $W$ and $V/W$. The latter forces $v \in W$ and the former $v = 0$, a contradiction. So $\fc = \fc_1$ or $\fc = \fc_2$.

Suppose $\fu = 0$. Then $\fc_2$ acts with a finite kernel on $W$, so $\rk \fc_2 \leq \rk W = 1$. Likewise $\rk \fc_1 \leq \rk (V/W) = 1$. Since $\fc$ equals $\fc_1$ or $\fc_2$, we have $\rk \fc \leq 1$ and $\rk \fg \leq 3$. By~\ref{t:module2:i:3}, $\fg \simeq \fsl_2(\bK)$ in the natural action. But there, $\fc$ \emph{is} a non-trivial $V$-nilpotent subring: a contradiction.
\item
Any $V$-nilpotent ring is easily seen to be $V$-nilpotent of class $\leq 2$, hence abelian.
If $V$ is $\fu$-irreducible, then $\fu$ acts freely on $V$, against $V$-nilpotence. So there is a $\fu$-series $0 < W < V$. The action of $\fu$ on $W$ is by (possibly trivial) scalars; , $V$-nilpotence forces triviality, so $W \leq C_V^\circ(\fu)$ and equality holds. In particular, $\fu$ may not normalise another non-trivial, proper submodule of $V$; which proves uniqueness. Always by $V$-nilpotence, $\fu$ centralises $V/W$. Thus $\fu \sqsubseteq C_\fg^\circ(W, V/W)$. The right-hand is $V$-nilpotent, so maximality of $\fu$ forces equality.
\item
Clearly $N_\fg^\circ(\fu) \sqsubseteq N_\fg^\circ(W_\fu)$. Moreover $N_\fg^\circ(W_\fu)$ is soluble by reducibility.
(For full rigour: apply~\ref{t:module2:i:soluble}. One may observe that $N_\fg^\circ(W_\fu)$ meets $C_\fg^\circ(w)$ at least in $\fu \neq 0$ for any $w \in W\setminus \{0\}$, or that this direction of ~\ref{t:module2:i:soluble} did not require meeting a centraliser.)
So there is a Borel subring $\fb \sqsubset \fg$ containing $N_\fg^\circ(W_\fu)$. Now $\fu \sqsubseteq N_\fg^\circ(\fu) \sqsubseteq \fb$, so $\fb$ meets a point centraliser. By~\ref{t:module2:i:soluble}, $\fb$ is reducible; there is a $\fb$-series $0 < W < V$. Then $\fu \sqsubseteq \fb$ acts on $W$, so $W = W_\fu$ by~\ref{t:module2:i:Wu}. Thus $\fb$ normalises $W_\fu$.
By Lie's theorem in the form of Corollary~\ref{c:solubleaction}, $\fb'$ acts trivially on $W_\fu$ and $V/W_\fu$, so by~\ref{t:module2:i:Wu} one has $\fb' \sqsubseteq C_\fg^\circ(W_\fu, V/W_\fu) = \fu$. Therefore $\fu$ is an ideal of $\fb$, viz.~$\fb \sqsubseteq N_\fg^\circ(\fu)$.
\qedhere
\end{enumerate}
\end{proofclaim}

\begin{step}[rigid, then disjoint, centralisers]\leavevmode
\begin{enumerate}[resume*=claims]
\item\label{t:module2:i:Cw}
For $w \in W_\fu \setminus\{0\}$, one has $C_\fg^\circ(w) = C_\fg^\circ(W_\fu)$.
\item\label{t:module2:i:u1u2}
There are two maximal distinct $V$-nilpotent subrings.
\item
We are done.
\end{enumerate}
\end{step}
\begin{proofclaim}\leavevmode
\begin{enumerate}[resume*=proofs]
\item
Let $w \in W_\fu\setminus \{0\}$ and $\fc = C_\fg^\circ(w)$ (not the same as in~\ref{t:module2:i:u}). By~\ref{t:module2:i:Cvsoluble}, $\fc$ is soluble, hence reducible by~\ref{t:module2:i:soluble}. Now $\fu \sqsubseteq \fc$, so every $\fc$-module is $\fu$-invariant. By~\ref{t:module2:i:Wu} the only possible non-trivial, proper $\fc$-submodule of $V$ is $W_\fu$. If $\fc$ acts on $W_\fu$ non-trivially, then it acts freely: a contradiction. So $\fc$ acts trivially on $W_\fu$, and $\fc \sqsubseteq C_\fg^\circ(W_\fu) \sqsubseteq C_\fg^\circ(w) = \fc$.
\item
Suppose not. Let $\fu$ be the unique such; let $W = C_V^\circ(\fu)$ and $\fb = N_\fg^\circ(\fu)$.
However intuitively pathological, this configuration takes a little time to dispose of.
We use~\ref{t:module2:i:Wu} and~\ref{t:module2:i:bu} implicitly.

Let $v_2 \in V\setminus C_V(\fu)$ and $\fc_2 = C_\fg^\circ(v_2)$, so that $\cork \fc_2 \leq 2$. By~\ref{t:module2:i:u}, $\fc_2$ contains a non-trivial, $V$-nilpotent subring $\fa_2 \neq 0$. By assumption, $\fa_2 \sqsubseteq \fu$. Now $\fc_2$ being soluble by~\ref{t:module2:i:Cvsoluble} must normalise a $1$-dimensional submodule by~\ref{t:module2:i:soluble}, say $W'$. Then $\fa_2$ centralises $W'$; if $W' \neq W$, then $\fa_2 \neq 0$ centralises $W + W' = V$, a contradiction. Hence $W' = W$ and $\fc_2 \leq N_\fg^\circ(W) = \fb$.

Since $v_2 \notin C_V(\fu)$, we have $v_2 \notin W$. So the action of $\fc_2$ on $V/W$ is not free; since $\rk(V/W) = 1$, the action is therefore trivial, viz.~$\fc_2 \leq C_\fg^\circ(V/W)$. Now $v_2 \notin C_V(\fu)$ implies $\fu \not\sqsubseteq \fc_2$, and therefore $\fc_2 \sqsubsetneq \fb \sqsubsetneq \fg$. Then $\rk \fc_2 = \rk \fg - 2$ and $\rk \fb = \rk \fg - 1$. It follows $\fb = \fu + \fc_2 \leq C_\fg^\circ(V/W) < \fg$ so $\fb = C_\fg^\circ(V/W)$.

Now recall that $\fu$ is $\fg$-nilpotent, viz.~$\ad$-nilpotent as a subring of $\fg$. Since $\rk(\fg/\fb) = 1$, the ring $\fu$ acts trivially on $\fg/\fb$.
Let $u \in \fu\setminus\{0\}$. Then $[u, \fg]\leq \fb \leq C_\fg(V/W)$.
Let $x \in \fg$, and $v \in V$.
Then:
\begin{align*}
xu \cdot v & = [x, u]\cdot v + u\cdot x v\\
& \quad \in [u, \fg] \cdot v + u \cdot V\\
& \qquad \leq C_\fg(V/W) \cdot v + W \leq W.\end{align*}
Since $0 < u V \leq W$ implies $u V = W$,
this yields $x W \leq W$. So any $x \in \fg$ normalises $W$, against irreducibility.
\item
By~\ref{t:module2:i:u1u2}, let $\fu_1 \neq \fu_2$ be two distinct maximal $V$-nilpotent subrings.
Let $W_i = C_V^\circ(\fu_i)$ be given by~\ref{t:module2:i:Wu}; notice $W_1 \neq W_2$.
Fix $w_i \in W_i\setminus\{0\}$ and let $\fc_i = C_\fg^\circ(w_i)$.
By~\ref{t:module2:i:Cw}, $\fc_1 \cap \fc_2$ centralises $W_1 + W_2 = V$, so $\fc_1 \cap \fc_2 = 0$. Since $\cork \fc_i = \rk (\fg\cdot v_i) \leq 2$, this gives:
\[2 \rk(\fg - 2) \leq \rk \fc_1 + \rk \fc_2 \leq  \fg,\]
whence $\rk \fg \leq 4$. By~\ref{t:module2:i:3} and~\ref{t:module2:i:4}, we are done.
\qedhere
\end{enumerate}
\end{proofclaim}
This completes the proof of Theorem~\ref{t:module2}.
\end{proof}

\begin{remark}[on \ref{t:module2:i:3}]
Presumably the following should hold by Chevalley's basis theorem.
Let $\fg$ be a connected Lie ring of finite Morley rank. Suppose $\fg$ is perfect and 
there is an algebraically closed field $\bK$ such that $\fg/Z(\fg)$ is a simple, finite-dimensional $\bK$-Lie algebra \emph{of Chevalley type} ($A_n, \dots, G_2$). If the characteristic is large enough, then $Z(\fg) = 0$.
This would be an analogue to \cite{ACCentral}. We do not know what happens with algebras `of Lie-Cartan type'.
\end{remark}

\subsection{Rewriting a classic}

We also include a streamlined proof of \cite[Theorem~A]{DActions}. The first author apologises for publishing, in his unipotent years, a less direct argument.

\begin{proof}
If $\chara V = 0$ then we linearise by \cite{DZilber}: there is a definable field $\bK$ such that $V \in \Vect{\bK}$ and $\fg \hookrightarrow \fgl(V\colon \Vect{\bK})$, definably.
If $\rk \bK = 2$ then the linear dimension is $1$ and we are done. If $\rk \bK = 1$, then $G= \GL_2(\bK)$ in its natural action, or $\rk G = 3$. In that case either $\fg = \SL_2(\bK)$ in its natural action, or $\rk (G \cap \SL_2(\bK))^\circ = 2$. But the latter is a normal subgroup of $G$, so $G$ is soluble and we are done.

We therefore suppose $\chara V > 0$ and denote it by $p$. Then we follow the proof of Theorem~\ref{t:module2}.

\ref{t:module2:i:rightkernel} is dealt with by \cite[Fact~2.7]{DActions}.

\ref{t:module2:i:3} does not use the recent \cite{FSimple}: indeed, if there is a centraliser with $\rk C_G^\circ(v) = 2$, then $G/Z(G) \simeq \PGL_2(\bK)$ by \cite{CGroups}. If all centralisers have $\rk C_G^\circ(v) = 1$, then all orbits $G\cdot v$ are generic, which creates an involution in $V$; so $G/Z(G) \simeq \PGL_2(\bK)$ by \cite{BBCInvolutions} and \cite{NNonsolvable}.

\ref{t:module2:i:4} does not use the recent \cite{WRank4}. If there is a centraliser with $\rk C_G^\circ(v) = 3$, then \cite{HAlmost} gives a normal subgroup and we linearise. If not, then $\rk C_G^\circ(v) = 2$ and we obtain an involution in $G$; so here again, pathological configurations are ruled out.

Now introduce the following.
\begin{definition*}
A definable, connected subgroup $U \leq G$ is \emph{$V$-unipotent} if for all $a \in U$, there is an integer $n$ such that $(a-1)^n = 0~[\End(V)]$.
\end{definition*}

In the present case, $n \leq 2 \leq p$, so $a^p - 1 = (a-1)^p = 0~[\End(V)]$ and $U$ has exponent $p$. Moreover, $U$ is abelian.

\ref{t:module2:i:u1u2} is trivial, as opposed to the Lie ring case.
\end{proof}

It would be interesting to translate \cite[Theorem~B]{DActions} for Lie rings. This would certainly use \cite[Variations~11 and~14]{DQuadratic}. The thorough discussion of pathological `partially quadratic' $\fsl_2(\bK)$-modules in characteristic $3$ \cite[\S~4.3]{DQuadratic} is not expected to have an impact here.

\paragraph{Acknowledgements.}
After \cite{DTLie} which was carried in Luminy and in Franceville thanks to \textsc{cirm}, we were given another opportunity to meet (this time in Paris) by \textsc{cirm} and the \textsc{ems}. Many thanks to both institutions; extra, personal thanks to Olivia Barbarroux, Pascal Hubert, and Carolin Pfaff.

We also greatly benefited from Ulla Karhumäki's energetic thoughts.
The first author completed his share of the work in the Mathematics Village of \c{S}irince, Turkey; best regards to Ali Nesin.

\printbibliography

\end{document}